\documentclass{amsart}
\usepackage{epsf}

\newtheorem{thm}{Theorem}[section]

\newtheorem{prop}[thm]{Proposition}

\theoremstyle{definition}
\newtheorem{defn}{Definition}

\theoremstyle{remark}
\newtheorem{rem}{Remark}

\newcommand{\R}{{\bf R}}
\newcommand{\C}{{\bf C}}
\newcommand{\Z}{{\bf Z}}
\newcommand{\half}{{\frac{1}{2}}}

\newcommand{\RPN}[1]{{\bf RP}^{#1}}
\newcommand{\CPN}[1]{{\bf CP}^{#1}}

\newcommand{\Rat}[2]{{\rm Rat}_{#1}({#2})}
\newcommand{\RRatI}[3]{{\rm RRat}_{#1,#2}({#3})}
\newcommand{\RRat}[2]{{\rm RRat}_{#1}({#2})}
\newcommand{\rat}[2]{{\bf rat}_{#1}({#2})}
\newcommand{\supp}[1]{{\rm supp}({#1})}

\begin{document}

\title[Spaces of rational loops]{Spaces of rational loops on a \\real 
projective space}

\author{Jacob Mostovoy}

\address{Instituto de Matem\' aticas (Unidad Cuernavaca), Universidad
Nacional Aut\' onoma de M\' exico,
Av. Universidad s/n, Col. Lomas de Chamilpa,
Cuernavaca, Morelos, M\' exico C.P. 62210}

.


%%%%%%%%%%%%%%%%%%%%%%%%%%%%%%%%%%%%%%%%%%%%%%%%%%%%%%%%%%%%%%%%%%%%%%%%%%
%%									%%
%%				ABSTRACT				%%
%%									%%
%%%%%%%%%%%%%%%%%%%%%%%%%%%%%%%%%%%%%%%%%%%%%%%%%%%%%%%%%%%%%%%%%%%%%%%%%%

\begin{abstract}
We show that the loop spaces on real  projective spaces are
topologically approximated by the spaces of rational maps
$\RPN{1} \rightarrow \RPN{n}$.  As a byproduct
of our  constructions we obtain an interpretation
of the Kronecker characteristic (degree) of an ornament via particle spaces.
\end{abstract}

\maketitle

%%%%%%%%%%%%%%%%%%%%%%%%%%%%%%%%%%%%%%%%%%%%%%%%%%%%%%%%%%%%%%%%%%%%%%%%%%
%%									%%
%%				INTRODUCTION				%%
%%									%%
%%%%%%%%%%%%%%%%%%%%%%%%%%%%%%%%%%%%%%%%%%%%%%%%%%%%%%%%%%%%%%%%%%%%%%%%%%

\section{Introduction.}
It is well-known that rational maps $\RPN{1}\rightarrow\RPN{m}$ are dense in
the space of all continuous maps from a circle to $\RPN{m}$.
Here we will show that the space of all such (basepoint-preserving) rational
maps is, in fact, homotopy equivalent to the loop space $\Omega\RPN{m}$.
For spaces of maps given by polynomials of bounded degree, we will
obtain an estimate of ``how well'' these spaces approximate $\Omega\RPN{m}$.
We consider maps which respect basepoints for convenience only;
corresponding statements for basepoint-free rational maps and free loops on 
$\RPN{m}$ easily follow.

Our results can be viewed as a real version of Segal's theorem \cite{Seg}.
Segal showed that the second loop space $\Omega^2\CPN{m}$
on a complex projective space can be topologically approximated
by the spaces of basepoint - preserving rational maps from $\CPN{1}$ to 
$\CPN{m}$. Namely, let $\Rat{n}{m}$ denote the space of
such rational maps of degree $n$ that send the point
$\infty\in\CPN{1}$ to some fixed point of $\CPN{m}$ with the topology of
a subset of $\Omega^2\CPN{m}$. It is convenient to choose
$(1,1,\ldots,1)$ as a basepoint in $\CPN{m}$; then every basepoint - 
preserving rational map can be given by a collection of monic polynomials
of the same degree. 
Let $(\Omega^2\CPN{m})_n$ be the
component of $\Omega^2\CPN{m}$ which parametrizes maps of degree $n$. Then 
the following is true:
\begin{thm}{\rm \cite{Seg}}\label{thm:c}
The natural inclusion
\[ \Rat{n}{m} \hookrightarrow (\Omega^2\CPN{m})_n \]
is a homotopy equivalence up to dimension $n(2m-1)$.
\end{thm}

This result was later generalized to spaces of maps into a wide class
of complex manifolds, see \cite{BHMM,G}.

Segal has also proved a real version of this theorem, which is different
from ours. Denote by
$\RRat{n}{m}$ the subspace of $\Rat{n}{m}$ of maps which commute with 
complex conjugation. Let $\RRatI{n}{k}{1}$ be the component of $\RRat{n}{1}$
which parametrizes maps whose restriction to the real line has degree $k$. 
Let also \linebreak $\Omega^2(\CPN{m},\RPN{m})$ be the
space of continuous basepoint - preserving maps
\[ (D^2, S^1, \ast) \rightarrow  (\CPN{m},\RPN{m}, \ast).\]
Such maps can be thought of as continuous conjugation-equivariant maps
\[\CPN{1}\rightarrow\CPN{m}.\]
Finally, let $\Omega^2(\CPN{1},\RPN{1})_{n,k}$ be the connected component
of $\Omega^2(\CPN{1},\RPN{1})$ which parametrizes maps that have degree $n$ 
as maps $\CPN{1}\rightarrow\CPN{1}$ and degree $k$ when restricted to the real
axis.
\begin{thm}{\rm \cite{Seg}}\label{thm:rc}
The natural inclusion
\[ \RRatI{n}{k}{1} \hookrightarrow \Omega^2(\CPN{1},\RPN{1})_{n,k} \]
is a homotopy equivalence up to dimension $\half(n-|k|)$.
\end{thm}

Now we will describe our real analogue of Theorem~\ref{thm:c}.

The restriction of any map  $f\in\RRat{n}{m}$ to the real axis defines a map
\[f_{\R}:\RPN{1}\rightarrow\RPN{m}.\]
In fact, $f_{\R}$ uniquely determines $f$,
so there is an inclusion $\RRat{n}{m}\hookrightarrow \Omega\RPN{m}$.
The image of this inclusion is not closed in $\Omega\RPN{m}$. For example, 
the sequence of maps
$g_{\epsilon}=\frac{x^2+1+\epsilon}{x^2+1}$ 
 does not have a limit in $\RRat{2}{1}$ as $\epsilon$ tends to 0.
(Here  $x$ is the coordinate on the affine piece of $\RPN{1}$.)
However, inside $\Omega\RPN{1}$ it converges to a constant map.

\begin{defn}
The space $\rat{n}{m}$ of real rational maps from $\RPN{1}$ to $\RPN{m}$ 
is the closure of $\RRat{n}{m}$ in $\Omega\RPN{m}$. 
\end{defn}

From what follows it will become clear that all points of $\rat{n}{m}$ 
correspond to rational
maps given by polynomials of degrees $k\leq n$, where $k$ and $n$ have the 
same parity. 

First we will look at the case $m=1$. 
\begin{thm}\label{thm:r1}
The space $\rat{n}{1}$ consists of $n+1$ contractible components,
indexed by the topological degree of the map $\RPN{1}\rightarrow\RPN{1}$ which
ranges from $-n$ to $n$ and has the same parity as $n$. 
\end{thm}
\begin{rem}
The components of $\Omega\RPN{1}$ are contractible.
\end{rem}

When $m>1$ the space  $\rat{n}{m}$ is connected. 
Denote by $(\Omega\RPN{m})_k$, where $k\in\Z/2$,
the connected components of $\Omega\RPN{m}$.
Then $\rat{n}{m}\subset (\Omega\RPN{m})_{n\ {\rm mod}\ 2}$.
\begin{thm}\label{thm:r}
For $m>1$ the natural inclusion
\[ \rat{n}{m} \hookrightarrow (\Omega\RPN{m})_{n\ {\rm mod}\ 2}\]
is a homotopy equivalence up to dimension $n(m-1)$.
\end{thm}
The proof goes along the lines of Segal's proof of Theorem~\ref{thm:c}.
In section 2 we define two types of configuration spaces which describe
divisors associated to real rational maps; we also define the stabilized 
spaces $\rat{\infty}{m}$.
Section 3 deals with maps from $\RPN{1}$ to $\RPN{1}$.
In section 4 we describe the stabilization results which prove
Theorem~\ref{thm:r} for $m>2$.
In section 5 we consider the action of $\pi_1 \rat{n}{2}$ on the higher
homotopy groups and finish the proof of Theorem~\ref{thm:r}.
Section 6 has no direct relevance to rational maps and is a byproduct
of the  constructions in section 5: we define an integer-valued invariant
of ornaments (see \cite{V}) which is interpreted as the degree of a
certain map $S^2\rightarrow S^2$.

%%%%%%%%%%%%%%%%%%%%%%%%%%%%%%%%%%%%%%%%%%%%%%%%%%%%%%%%%%%%%%%%%%%%%%%%%%
%%									%%
%%		CONFIGURATIONS OF REAL DIVISORS				%%
%%									%%
%%%%%%%%%%%%%%%%%%%%%%%%%%%%%%%%%%%%%%%%%%%%%%%%%%%%%%%%%%%%%%%%%%%%%%%%%%

\section{Configurations of real divisors.} 
Maps in $\RRat{n}{m}$ can be parametrized by collections of $(m+1)$ 
positive real divisors of degree $n$ in \C, subject to the following 
condition: the intersection of the supports of all $m+1$ divisors is empty. 

If we wish to find a configuration space which would parametrize maps in
$\rat{n}{m}$, we must
substitute this condition for a weaker one, namely: the intersection of
supports of all $m+1$ divisors  must be disjoint from the real axis.
We denote the space of such collections of divisors by $D(n,m)$.

Choose $z\in\C,\ Im(z)\neq 0$. There is an inclusion map 
\[ D(n,m) \stackrel{i'_n}{\hookrightarrow} D(n+2,m), \]
defined by adding the divisor $z+\bar{z}$ to each of the $m+1$ divisors.

Notice that we do not obtain a 1-1 correspondence between $\rat{n}{m}$
and $D(n,m)$; for example, pairs of polynomials $(x^2+1,x^2+1)$ and 
$(x^2+2,x^2+2)$ define the same (constant) map in $\rat{2}{1}$.
However, we have a continuous map
\[ P: D(n,m) \rightarrow  \rat{n}{m}.\]
Roughly speaking, if we think of points of $D(n,m)$ as of collections of 
particles of several kinds, we allow particles of different kinds to cancel out
provided they are not on the real axis.
The effect of this cancellation is that the degree of all divisors
drops by 2, so we have a commutative diagram
\[ \begin{array}{ccc}
	D(n,m)		&\stackrel{i'_n}{\hookrightarrow} 	& D(n+2,m)\\
	\downarrow P	&				& \downarrow P\\
	\rat{n}{m} 	&\stackrel{i_n}{\hookrightarrow}	& \rat{n+2}{m}
\end{array} \]
Notice that the inclusion map $i_n$ is canonical.
Clearly, 
\[ \rat{n}{m}\backslash\rat{n-2}{m} = \RRat{n}{m}.\]

There is another space of real divisors which will be of interest for us.
Let $\R^{(n)}$ denote the $n$-fold symmetric product of $\R$.
We write points in $\R^{(n)}$, i.e. divisors on $\R$
as linear combinations
$\sum k_ix_i$ of points $x_i\in \R$ with non-negative (but possibly zero)
multiplicities. (Here, obviously, $\sum k_i=n$.)
Define the subset $S(n,m)$ of the $(m+1)$-fold Cartesian product
$[\R^{(n)}]^{m+1}$  by the condition that the intersection of
supports of all $m+1$ divisors is empty.
We define the configuration space $T(n,m)$ of mod 2 particles on
the real line
as the quotient of $S(n,m)$ by the following equivalence relation:
\[ (\sum_i k^{0}_{i}x^{0}_{i},\ldots,\sum_i k^{m}_{i}x^{m}_{i})\sim
(\sum_i l^{0}_{i}x^{0}_{i},\ldots,\sum_i l^{m}_{i}x^{m}_{i})\]
if  $k^{p}_{q}\equiv l^{p}_{q} \pmod{2}$ for all $p,q$.

The points of $T(n,m)$ are, essentially, collections of particles of
$(m+1)$ different kinds on \R,  such that two particles of the same kind 
annihilate each other and such that $m+1$ particles of different kinds
cannot collide in the same point.

Define a map 
\[ D(n,m)\rightarrow T(n,m) \]
by sending a collection of polynomials which is defined by a point in $D(n,m)$
to the collection of their real root systems with multiplicities 
reduced modulo 2.
It is easy to check that this map is continuous and factors through 
$\rat{n}{m}$.
There is a commutative diagram 
\[ \begin{array}{ccc}
\rat{n}{m} 	&\stackrel{i_n}{\hookrightarrow}	& \rat{n+2}{m}	\\
\downarrow Q	&					& \downarrow Q	\\
T(n,m)		&\stackrel{i''_n}{\hookrightarrow} 	& T(n+2,m)	
\end{array} \]
Both inclusions in this diagram are canonical.

There is another kind of inclusion maps 
\[ j'_n: D(n,m)\hookrightarrow D(n+1,m), \]	
defined by ``adding real points at infinity''. We can think of $D(n,m)$ 
as consisting of collections of divisors that are
contained in an open disk of radius $n$ with the center at $0$.
Choose  $x_i\in\R$, $0\leq i \leq m$ such that $n<x_0<x_1<\ldots<x_m<n+1$.
Then the inclusion $j'_n$ is defined by adding $x_i$ to the $i$-th
divisor in the collection for each $i$.

There are induced inclusions 
\[ j_n: \rat{n}{m}\hookrightarrow \rat{n+1}{m} \] and 
\[ j''_n: T(n,m)\hookrightarrow T(n+1,m).\]
None of these maps are canonical; it is easy to see, however, that
for $m>1$
$j'_{n+1}j'_n$ is homotopic to $i_n$, similarly for $j_n$ and $j''_n$.

Now for $m>1$ we can define the space $\rat{\mathrm{odd}}{m}$
as the direct limit of the sequence
\[ \rat{1}{m}\stackrel{i_1}{\hookrightarrow}\rat{3}{m}\stackrel{i_3}{\hookrightarrow}\ldots \]
The spaces $\rat{\mathrm{ev}}{m}$, $D(odd,m)$, $D(ev,m)$, $T(odd,m)$ and 
$T(ev,m)$ are defined similarly as direct limits.
By $\rat{\infty}{m}$ and  $T(\infty,m)$ we denote the disjoint unions 
$\rat{\mathrm{odd}}{m}\cup\rat{\mathrm{ev}}{m}$
and $T(odd,m)\cup T(ev,m)$ 
respectively.  As before, the collections of divisors in $D(ev, m)$ can be 
thought of as lying in an open disk of radius 1 with the center at 0. 
Adding a point $x_k=2+k$ to the $k$th divisor for all $0\leq k\leq m$ we get
a map $j_{\mathrm{ev}}'$ from $D(ev,m)$ to  $D(odd,m)$, which is clearly
a homotopy equivalence. This map also descends to homotopy equivalences
\[ j_{\mathrm{ev}}:\rat{\mathrm{ev}}{m}\rightarrow\rat{\mathrm{odd}}{m} \]
and 
\[ j_{\mathrm{ev}}'':T(ev, m)\rightarrow T(odd, m). \]

The spaces $D(n,m)$, $\rat{n}{m}$ and $T(n,m)$ have homotopy
types of CW-complexes. This is clear for $D(n,m)$ as it is a smooth manifold; 
for $\rat{n}{m}$ and $T(n,m)$ it can be established by an argument 
similar to the proof of Lemma~2.2 of \cite{McDuff}.

The reason for considering the configuration spaces $D(n,m)$ and $T(n,m)$
is the following
\begin{prop}\label{prop:main}
The maps 
\[ D(n,m)\stackrel{P}{\longrightarrow}\rat{n}{m} \stackrel{Q}{\longrightarrow}
 T(n,m) \]
are homotopy equivalences for all $n$ and $m$, including $n=\infty$.
\end{prop}
\begin{proof}
Essentially, the reason $P$ and $Q$ are homotopy equivalences is that they
have contractible fibers.
As we have seen above, the space $\rat{n}{m}$ has a canonical filtration
by subspaces $\rat{k}{m}$, where $k\leq n$ and $k\equiv n\pmod{2}$.
The fiber of $P$ over $\rat{k}{m}\backslash\rat{k-2}{m}=\RRat{k}{m}$
is the $\frac{1}{2}(n-k)$-fold symmetric product of the open upper half-plane,
i.e. $\R^{n-k}$. In fact, $P^{-1}(\RRat{k}{m})$ is the direct
product $\RRat{k}{m}\times\R^{n-k}$. The space $D(n,m)$ can be triangulated
in such a way that the inverse image $P^{-1}(\rat{k}{m})$ is 
a subcomplex, so the inclusion $P^{-1}(\rat{k}{m})\hookrightarrow D(n,m)$ 
is a cofibration.  Let $p_k$ be the sequence of maps which collapse
the fibers of $P$ over $\RRat{k}{m}$ consecutively. It is a sequence of 
homotopy equivalences. On the other hand, $P=p_{n}p_{n-2}\ldots$

The argument for $T(n,m)$ is similar. The space $T(n,m)$ is graded by
subspaces $T_{\bf k}(n,m)$, where ${\bf k}=(k_0,\ldots,k_m)$ is a set
of integers such that $0\leq k_i\leq n$ and $k\equiv n\pmod{2}$.
A point of $T(n,m)$ belongs to $T_{\bf k}(n,m)$ if it has exactly
$k_i$ particles of the $i$th kind.

The closure of $T_{\bf k}(n,m)$ is contained in the union
$\cup_{{\bf l}\leq{\bf k}}T_{\bf l}(n,m)$; here by ${\bf l}\leq{\bf k}$
we mean $l_i\leq k_i$ for any $i$.
The restriction of $QP$ to $P^{-1}Q^{-1}(T_{\bf k}(n,m))$ is clearly
a homotopy equivalence and, for the same reason as before, the inclusion
\[ P^{-1}Q^{-1}(\cup_{{\bf l}\leq{\bf k}}T_{\bf l}(n,m))
\hookrightarrow D(n,m) \]
is a cofibration.
Collapsing the fibers of $QP$ consecutively we get a chain of homotopy
equivalences so $QP$ (and, hence, $Q$) is itself a homotopy equivalence.
\end{proof}
As an immediate corollary we get that
for $m>1$ spaces $\rat{n}{m}$ are connected and for $m>2$ - simply-connected.
Indeed, $D(n,m)$ is obtained by removing a triangulable subset of 
codimension $m$ from $\R^{n(m+1)}$.

%%%%%%%%%%%%%%%%%%%%%%%%%%%%%%%%%%%%%%%%%%%%%%%%%%%%%%%%%%%%%%%%%%%%%%%%%%
%%									%%
%%		SPACES OF MAPS \RPN{1} \rightarrow \RPN{1}		%%
%%									%%
%%%%%%%%%%%%%%%%%%%%%%%%%%%%%%%%%%%%%%%%%%%%%%%%%%%%%%%%%%%%%%%%%%%%%%%%%%

\section{Spaces of maps $\RPN{1} \rightarrow \RPN{1}$.} 

It is well-known that $\RRat{n}{1}$ has $n+1$ connected components.
They are indexed by the degree of maps $\RPN{1}\rightarrow\RPN{1}$, which
ranges from $-n$ to $n$ and has the same parity as $n$, see \cite{Bro} and
\cite{Seg}. The same is true for $\rat{n}{1}$ as  
$\RRat{n}{1}\subset\rat{n}{1}$ is a dense subset.
Alternatively, it is easy to see directly that $T(n,1)$ has exactly $n+1$ 
components.

Theorem~\ref{thm:r1} follows from
\begin{prop}\label{prop:ratn1}
Connected components of $T(n,1)$ are contractible.
\end{prop}
\begin{proof}
Points of $T(n,1)$ are pairs of disjoint mod 2 divisors on a real line.
We can treat them as sets of particles of different ``charges'', such
that particles of the same charge annihilate each other and particles of
different charges cannot collide. Let us extend the analogy.

At the time $t=0$ let all the particles be at rest.
We introduce three forces: gravity with the potential $x^2$,
friction which is proportional to the velocity of the particle,
and interaction between particles. A particle interacts only with its
immediate neighbors; it attracts particles of the same charge and repels
particles of the opposite charge; the magnitude of this force is $\frac{1}{r}$,
where $r$ is the distance between particles. The dynamics of this system
defines a flow on $T(n,1)$. As particles of the same charge annihilate each
other, the interaction forces change discontinuously, but, as the force
is essentially the second derivative of the coordinate, the flow is
continuous, at least on any compact subset. Calculation shows that there is 
exactly one state of equilibrium for each connected component of $T(n,1)$, 
so as $t\rightarrow\infty$ we get a homotopy of any map $S^k\rightarrow T(n,1)$
to a constant map. Hence, all components of $T(n,1)$ have trivial homotopy 
groups and this implies that they are contractible.
\end{proof}

%%%%%%%%%%%%%%%%%%%%%%%%%%%%%%%%%%%%%%%%%%%%%%%%%%%%%%%%%%%%%%%%%%%%%%%%%%
%%									%%
%%  THE HOMOTOPY EQUIVALENCE OF \rat{\infty}{m} AND	\Omega\RPN{m} 	%%
%%									%%
%%%%%%%%%%%%%%%%%%%%%%%%%%%%%%%%%%%%%%%%%%%%%%%%%%%%%%%%%%%%%%%%%%%%%%%%%%

\section{The homotopy equivalence $\rat{\infty}{m}\hookrightarrow\Omega\RPN{m}$
and homology stabilization.}

\begin{prop}\label{prop:limit}
For all $m$ the natural inclusion
\[\rat{\infty}{m}\hookrightarrow\Omega\RPN{m}\]
is a homotopy equivalence.
\end{prop}

\begin{proof}
The statement of the Proposition for $m=1$ follows from Theorem~\ref{thm:r1},
so we will assume that $m\geq 2$.

Let $y_0,\ldots,y_m$ be the homogeneous coordinates
in $\RPN{m}$. Denote by $\Sigma$ the subspace of loops which have
an infinite number of intersection points with some of the hyperplanes
$y_k=0$. The inclusion map 
$\Omega\RPN{m}\backslash \Sigma \hookrightarrow \Omega\RPN{m}$ 
is a homotopy equivalence.
The reason for this is that $\Omega\RPN{m}$ can be thought of as an
infinite-dimensional manifold and $\Sigma$ is a subset of infinite
codimension, so removing $\Sigma$ doesn't change the homotopy type of
$\Omega\RPN{m}$. However, rather than making this argument precise,
we will sketch a direct proof.

Every loop on $\RPN{m}$ can be lifted to a loop or a path on
$S^{m}=\{y_0,\ldots,y_m |$ $y_0^2+\ldots+y_m^2=1\}.$
First with the help of some self-homeomorphism of $S^{m}$
we deform the space $\Omega\RPN{m}\backslash \Sigma$ into the space of loops 
(paths)
on $S^{m}$ that have not more than a finite number of intersections with 
any of the sets $y_k=\epsilon$ for some small $\epsilon$. 
Notice that any piecewise geodesic (in the standard metric) path on $S^m$ has 
only a finite number of intersections with $y_k=\epsilon$ for any $k$.
Now one can use piecewise geodesic approximations and the methods of 
\cite{Mil} can be applied. The argument is identical to
the proof of the Theorem~17.1 of \cite{Mil} with the space of piecewise 
smooth loops replaced by $\Omega\RPN{m}\backslash \Sigma$.

Define a map 
$q:\Omega\RPN{m}\backslash \Sigma \stackrel{q}{\longrightarrow} T(\infty,m)$ by
sending a loop $\omega$ to the collection of mod 2 divisors
$(\zeta_0,\ldots,\zeta_m)$,
where $\zeta_k$ consists of the values of the parameter $x$ at which
$\omega$ intersects the hyperplane $y_k=0$; the coefficient of  $x$
in $\zeta_k$ is $1$ if $\omega$ crosses $y_k=0$ near $x$ and $0$ otherwise.

In the commutative diagram
\[ \begin{array}{rrl}
\rat{\infty}{m} &\hookrightarrow 
&\Omega\RPN{m}\backslash \Sigma\ \simeq\ \Omega\RPN{m}	\\
\downarrow Q	&  \swarrow q  &\\
T(\infty,m)		&	&
\end{array} \]
$Q$ is a homotopy equivalence, 
so all maps 
$\pi_*\rat{\infty}{m}\hookrightarrow\pi_*\Omega\RPN{m}$ are injective.

On the other hand, the Weierstra\ss\  approximation theorem implies that
the above maps are surjective. Indeed, consider the inclusion of the even
components $\rat{\mathrm{ev}}{m}\hookrightarrow (\Omega\RPN{m})_0=\Omega S^m$. 
Each loop on $S^m$
can be represented by $m+1$ functions $f_0(x),\ldots ,f_m(x)$ such that
$\lim_{x\rightarrow\infty}{f_i(x)}=1$ for any $i$ and $\sum_0^m f_i(x)^2 =m+1$.
If we identify $\R\cup\{\infty\}$ with $S^1$ by substituting
$x=\tan{\frac{\alpha}{2}}$ the functions $f_i$ become continuous functions
on $S^1$ and, hence, can be uniformly approximated by trigonometric 
polynomials.
It is easy to check that trigonometric polynomials in $\alpha$ are rational
functions in $x$, so any loop on $S^m$ can be uniformly approximated by
rational loops.
The same argument establishes that any map $S^k\rightarrow(\Omega\RPN{m})_0$ 
can be
uniformly approximated by maps $S^k \rightarrow \rat{\mathrm{ev}}{m}$. This is 
enough to
claim that every element of $\pi_k(\Omega\RPN{m})_0$ is an image of some 
element of $\pi_k\rat{\mathrm{ev}}{m}$. The same is true for the odd 
components as there is a homotopy commutative diagram
\[ \begin{array}{ccc}
\rat{\mathrm{ev}}{m} &\stackrel{j_{\mathrm{ev}}}{\rightarrow} 
&\rat{\mathrm{odd}}{m}\\
\downarrow 	&    			&\downarrow \\
(\Omega\RPN{m})_0 &\simeq 	&(\Omega\RPN{m})_1.	
\end{array} \]
\end{proof}
\begin{prop}\label{prop:stability}
The inclusion map 
$j_n: \rat{n}{m}\hookrightarrow\rat{n+1}{m}$ is a homology
equivalence up to dimension $n(m-1)$
for all $m>1$ and all $n>0$. 
\end{prop}
As $j_{n+1}j_n$ is homotopic to the inclusion 
$i_n:\rat{n}{m}\hookrightarrow\rat{n+2}{m}$
it is enough to prove that $i_n$ is a homology equivalence up to dimension 
$n(m-1)$. In view of Proposition~\ref{prop:main} this follows from the 
corresponding statement for the map $i'_n: D(n,m)\hookrightarrow D(n+2,m)$. 
The proof of the latter repeats Segal's proof of
Proposition~5.1 in \cite{Seg} almost  word-for-word; we refer to \cite{Seg}
for details.

In the case when $m>2$ $\rat{n}{m}$ is simply-connected and so $j_n$
is a homotopy equivalence up to dimension $n(m-1)$. In order to
show that this is also the case for $m=2$ it is necessary to check that
$\rat{n}{2}$ is simple up to dimension $n$, i.e. that for $k<n$
the action of $\pi_1\rat{n}{2}$ on $\pi_k\rat{n}{2}$ is trivial.

%%%%%%%%%%%%%%%%%%%%%%%%%%%%%%%%%%%%%%%%%%%%%%%%%%%%%%%%%%%%%%%%%%%%%%%%%%
%%									%%
%%		SPACES OF MAPS $\RPN{1}\rightarrow\RPN{2}$.		%%
%%									%%
%%%%%%%%%%%%%%%%%%%%%%%%%%%%%%%%%%%%%%%%%%%%%%%%%%%%%%%%%%%%%%%%%%%%%%%%%%

\section{Spaces of maps $\RPN{1}\rightarrow\RPN{2}$.}
First let us consider the spaces $\rat{1}{2}$ and $\rat{2}{2}$.
The space $D(1,2)$ is $\R^3$ with the diagonal line 
 removed, so $\rat{1}{2}\simeq S^1$. A generator of
$\pi_1 D(1,2)$ can be represented by a loop 
\[ t \rightarrow (x+\cos{t}, x+\sin{t}, x);\ t\in [0,2\pi]. \]
 
The space $\rat{2}{2}$ is homotopy equivalent to a torus with one of the 
meridians collapsed so, in particular, $\pi_1\rat{2}{2}=\Z$. This can be
seen as follows.
 
Points in $D(2,2)$ can be treated as maps $\R\rightarrow\R^3$. Images of these
maps are parabolas (possibly degenerated to rays) which miss 0 and whose
asymptotic direction is $(1,1,1)$. Let $y_1, y_2, y_3$ be the coordinates in 
$\R^3$. The map from $D(2,2)$ to 
$(\R^2\backslash 0)\times(\R^2\backslash 0)$ with the diagonal collapsed
is given by sending the first point of intersection of the parabola
with the plane $y_1+y_2+y_3=0$ to the first copy of $\R^2\backslash 0$,
and the second point of intersection to the second copy. If the parabola
does not intersect the plane at all, we map it into the collapsed diagonal.
This map is a quasifibration with contractible fibres and, hence
a homotopy equivalence.

A generator of
$\pi_1 D(2,2)$ can be represented by a loop 
\[ t \rightarrow 
(x^2+x+\cos{t}+\frac{1}{2}\sin{t},\  x^2+x+\cos{t}-\frac{1}{2}\sin{t},
\ x^2-x+\cos{t}), \]
where $t\in [0,2\pi]$.

\begin{prop}\label{prop:pi1}
$\pi_1\rat{n}{2}=\Z$ for any $n>0$. 
\end{prop}
\begin{proof}
The loops of even degree in $\RPN{2}$ can be lifted to  loops on a 2-sphere; 
(we think of the basepoint as of the  north pole); loops of odd degree can be 
lifted to paths from the north pole to the south pole. Fix a group structure 
on the equator $S^1$.
We define a map $\rat{n}{2}\rightarrow S^{1}$ by sending a rational map $f$ to
the alternating sum $b_1-b_2+b_3-\ldots$, where $b_i$ is the $i$-th point
of intersection of the image of $f$ with the equator. (Here we count
the intersection points with multiplicities so, for example, a double
point of intersection will have no effect on the sum.) It can be verified 
directly that the composite maps
\[ \pi_1\rat{1}{2}\stackrel{I_{2k+1}}{\longrightarrow}
\pi_1\rat{2k+1}{2}\rightarrow \pi_1 S^{1} \]
and
\[ \pi_1\rat{2}{2}\stackrel{I_{2k}}{\longrightarrow}\pi_1\rat{2k}{2}
\rightarrow\pi_1 S^{1} \]
are isomorphisms. Hence  the maps
$I_{2k+1}:\pi_1\rat{1}{2}\rightarrow\pi_1\rat{2k+1}{2}$
and
$I_{2k}:\pi_1\rat{2}{2}\rightarrow\pi_1\rat{2k}{2}$
are injective.

Notice that the equator considered as a group acts freely on the spaces
$\rat{2k+1}{2}$ by rotating the sphere. This action commutes with the map
$\rat{2k+1}{2}\rightarrow S^{1}$ defined above; moreover, the orbits map 
onto $S^1$ homeomorphically, so $\rat{2k+1}{2}$ splits as a
product of a circle with the quotient space by the action. This will
imply that the spaces $\rat{2k+1}{2}$ are homotopy simple in all dimensions
as soon as we have 
proved that $\pi_1 \rat{2k+1}{2}=\Z$. 
This argument fails for the loops of
even degree; and, indeed, one can show that $\pi_1\rat{2}{2}$ acts 
non-trivially on $\pi_2\rat{2}{2}$.
\begin{rem}
It is well-known that there
is a decomposition
\[ \Omega S^2 \simeq S^1\times \Omega S^3 \]
associated to the the Hopf fibration.
\end{rem}
Now, rather than looking at loops on $\rat{n}{2}$ we will
consider loops on $T(n,2)$. They can be presented in
a braid-like fashion; Figure~1\ shows  generators of 
$\pi_1 T(1,2)$ and $\pi_1 T(2,2)$. (A loop is obtained by moving the 
horizontal line down and registering the intersection with the ``braid'').

\begin{figure}
\[\epsffile{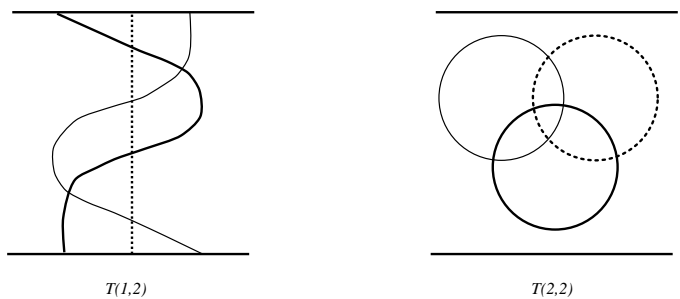}\]
\caption{Generators of $\pi_1 T(1,2)$ and $\pi_1 T(2,2)$.}
\label{fig:pi1}
\end{figure}

To prove that the maps $I_n$ are surjective we show directly that
any loop in $\rat{2k+1}{2}$ or $\rat{2k}{2}$ can be deformed into 
$\rat{1}{2}$ or $\rat{2}{2}$ respectively.

Let ${\zeta}(t)=(\zeta_1(t),\zeta_2(t),\zeta_3(t))$ be a loop on
$T(2k+1,2)$. 
By ``support of a family of mod 2 divisors $\zeta_i(t)$'', $i=1,2,3$,
we will mean the closure of the set $\supp{\zeta_i(t)}$ in $\R\times[0,1]$.
Without loss of generality assume that supports of
families of divisors $\zeta_2(t)$ and $\zeta_3(t)$ intersect
in a finite number of points $a_i$, corresponding to distinct values
of parameter $t_i$, where $0\leq i\leq C$ for some integer $C$.

First of all we deform ${\zeta}(t)$ in such a way that $\supp{\zeta_1(t)}$
for each fixed $t$ consists of one point. This is
done in two steps.

First we find such an $\epsilon$  that:

a) the support of the family
$\zeta_1(t)$ is disjoint from the rectangle
\[R=[a_i-\epsilon, a_i+\epsilon]\times[t_i-2\epsilon, t_i+2\epsilon]\] in
$\R\times [0,1]$ for any $0\leq i\leq C$;

b) $t_{i+1}-t_i>4\epsilon$ for $0\leq i < C$ and $t_0 > 2\epsilon$ and
$t_C < 1-2\epsilon$.

Let us fix $i$ for the moment.
 On the interval of the parameter $|t-t_i|< 2\epsilon$
the divisor $\zeta_1(t)$ can be expressed as a sum
$\zeta'_1(t)+\zeta''_1(t)$ such that $\zeta'_1(t)$ lies to the left
 and  $\zeta''_1(t)$ lies to the right
of  $a_i$. As the mod 2 degree of the divisor $\zeta_1(t)$ is 1 for any $t$,
either $\zeta'_1(t)$ or $\zeta''_1(t)$ has mod 2 degree 0 for any 
$t\in [t_i-2\epsilon, t_i+2\epsilon]$. Let, for example, the mod 2 degree
of $\zeta''_1(t)$ be equal to 0.

Let $N$ be so big that the support of the family $\zeta_1(t)$ lies
in $(-\infty, N)$.
Define a homotopy $\phi_u$ from $\R\times [0,1]\backslash R$ to
itself such that:

a) it is constant everywhere apart from the rectangle
$[a_i+\epsilon, N]\times [t_i-2\epsilon, t_i+2\epsilon]$;

b) it carries the rectangle 
$[a_i+\epsilon, N]\times [t_i-\epsilon, t_i+\epsilon]$ onto the interval
$N\times [t_i-\epsilon, t_i+\epsilon]$; see Figure~2.

\begin{figure}
\epsffile{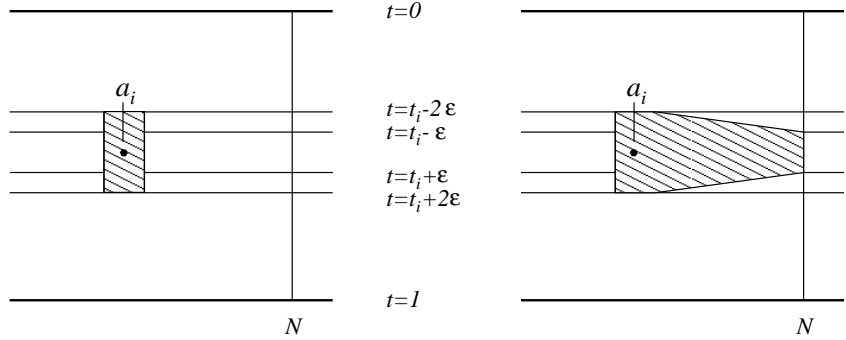}
\caption{The homotopy $\phi_u$ at $u=0$ and $u=1$.}
\label{fig:push}
\end{figure}

We change $\zeta(t)$ by deforming $\zeta_1(t)$ with the help of the homotopy
$\phi_u$ and leaving $\zeta_2(t)$ and $\zeta_3(t)$ fixed. The effect of this
operation is that the points of $\zeta_1(t)$ to the right of $a_i$ cancel out.

After performing this deformation for each $i$ we do the second step.
Let $g(t)$, $t\in[0,1]$ be a function such that:

a) for each $i$ $g(t_i)< a_i$ if the mod 2 degree of the divisors
$\zeta'_1(t)$ is 1, and $g(t_i) >  a_i$ otherwise;

b) the graph of $g(t)$   intersects the  support of
the family of divisors $\zeta_3(t_i)$
in a finite number of points. 

Now, let $\psi_u$ be a retraction  of $\R\times [0,1]$ onto the graph of $g(t)$
which preserves lines $t={\rm const.}$
We deform $\zeta(t)$ by pushing $\zeta_1(t)$ towards the curve $g(t)$
with the help of  $\psi_u$. 

The result of the above tricks is that for any $t$ the deformed  $\zeta_1(t)$
consists of one point only, and $\zeta_2(t)$ and $\zeta_3(t)$ don't
change. So after performing  similar operations with $\zeta_2(t)$ and 
$\zeta_3(t)$ we get a loop in $T(1,2)$ which is homotopic to $\zeta(t)$.

The proof that any loop in $T(2k,2)$ can be deformed into $T(2,2)$ is
similar, but easier; we leave it to the reader.
\end{proof}

We have already seen that the spaces $\rat{n}{2}$ are homotopy simple for
$n$ odd.
The proof of Theorem~\ref{thm:r} will be finished as soon as we have proved
\begin{prop}
For $n$ even the space  $\rat{n}{2}$ is homotopy simple up to dimension
$n$. 
\end{prop}
\begin{proof}
By Proposition~\ref{prop:pi1} $\pi_1\rat{n}{2}=\Z,$ so $\rat{n}{2}$ is homotopy
simple up to dimension 1.

Suppose we know that for some $k<n$ the space $\rat{n}{2}$
is homotopy simple up to dimension $k$. 
By Proposition~\ref{prop:stability} the map 
$j_{n-1}:\rat{n-1}{2}\hookrightarrow\rat{n}{2}$ is a homology equivalence up 
to dimension
$n-1$. This implies that $j_{n-1}$ is a homotopy equivalence up to dimension
$k$ and, in particular, that the map  
$\pi_k\rat{n-1}{2}\hookrightarrow\pi_k\rat{n}{2}$ 
is surjective. But, as $j_{n-1}$ induces an isomorphism of fundamental groups,
this means that the action of  $\pi_1\rat{n}{2}$ on $\pi_k\rat{n}{2}$
is trivial, and, hence, $\rat{n}{2}$
is simple up to dimension $k+1$. 
\end{proof}

%%%%%%%%%%%%%%%%%%%%%%%%%%%%%%%%%%%%%%%%%%%%%%%%%%%%%%%%%%%%%%%%%%%%%%%%%%
%%									%%
%%		THE DEGREE OF AN ORNAMENT				%%
%%									%%
%%%%%%%%%%%%%%%%%%%%%%%%%%%%%%%%%%%%%%%%%%%%%%%%%%%%%%%%%%%%%%%%%%%%%%%%%%
\section{The degree of an ornament.}	

A 3-ornament is a smooth map of the disjoint union of
3 circles into $\R^2$,  such that images of all three circles do not
pass through the same point, see \cite{V}. Two  3-ornaments are equivalent
if they are connected by a path in the space of such maps.

The most basic example of an (unoriented) ornament is shown on the 
right-hand side of Figure~1; we considered it as a generator of
$\pi_1 T(2,2)$. In fact, any 3-ornament defines an element of
$\pi_1 T(ev,2)$ in a similar fashion unless it has an infinite number
of intersections with some horizontal line. However, this is a condition
of infinite codimension so we can always deform a 3-ornament so that
it does define an element of $\pi_1 T(ev,2)=\Z$, which we call the
degree of the ornament. Clearly, equivalent ornaments have the same degree.

The reason for calling this invariant ``degree'' is the following
interpretation of it. Consider a regular ornament, i.e. such that it is
an immersion and all the multiple points of its image in $\R^2$ are
double transversal intersection points only. Let $f_1, f_2$ and $f_3$ be 
functions from $\R^2$ to $\R$, such that the zero set of each function defines
the corresponding component of the ornament (up to parametrization, of course) 
and the zero set of each of the functions contains only
a finite number of critical points. Assume also that at infinity 
$f_i\rightarrow 1$, $i=1,2,3$. Then we have a continuous map
\[ S^2 =\R^2\cup\{\infty\}\stackrel{f_1,f_2,f_3}
{\longrightarrow}\R^3\backslash 0 \simeq S^2.\]
The degree of the ornament is just the degree of this map.

Even though the above construction cannot be found in \cite{V}, historically
it precedes not only ornaments, but, in fact, the very notion of degree.
It can be found, in a more general setting, in Kronecker's paper \cite{Kro}, 
where he defines the "characteristic" of a system of functions with no common
zeros. Kronecker's definition of the characteristic of the system 
$(f_1,f_2,f_3)$ was roughly as follows.

Suppose that the zero sets of functions $f_i$ intersect transversally.
To each point of intersection $P_k$ of the curves $f_1=0$ and $f_2=0$ we assign
an integer which is equal to 0 if $f_3(P_k)<0$, and $\pm 1$ if $f_3(P_k)>0$
and the tangent vectors to the curves $f_1=0$ and $f_2=0$ at $P_k$ define a 
positive/negative
orientation of $\R^2$. (We assume that the orientations of the curves $f_i=0$
are induced by the standard orientation of $\R^2$.)
The characteristic of the system of functions $(f_1,f_2,f_3)$ is the sum of
these integers over all $P_k$.

It is clear  that this definition amounts to calculating the algebraic number 
of inverse images of a point in $S^2$, so the Kronecker
characteristic of a system $(f_1,f_2,f_3)$ is the same thing as the degree
of a map given by these functions. 
Modern references for the Kronecker characteristic are \cite{Hirsch} and,
less recent, but probably more precise, \cite{Lef}.

I would like to thank A.B. Merkov for remarking  that there is a 
characteristic-type formula for the degree of an ornament and V.I. Arnold who 
pointed out that it was, in fact, first discovered by Kronecker.

%%%%%%%%%%%%%%%%%%%%%%%%%%%%%%%%%%%%%%%%%%%%%%%%%%%%%%%%%%%%%%%%%%%%%%%%%%
%%									%%
%%		 BIBLIOGRAPHY						%%
%%									%%
%%%%%%%%%%%%%%%%%%%%%%%%%%%%%%%%%%%%%%%%%%%%%%%%%%%%%%%%%%%%%%%%%%%%%%%%%%

\end{document}